\documentclass[11pt,twoside]{article}

\usepackage{amsmath}
\usepackage{amsfonts}
\usepackage{amssymb}
\newcommand{\RR}{\mathbb R}

\newcommand{\email}[1]{{\small E-mail: {\textsf {#1}}}}
\newtheorem{defi}{Definition}

\newtheorem{theo}[defi]{Theorem}
\newtheorem{Le}[defi]{Lemma}
\newtheorem{Cor}[defi]{Corollary}
\newtheorem{Rem}[defi]{Remark}

\newcommand\qed{\hfill $\square$}

\begin{document}

\title{Global existence vs. blowup for the one dimensional quasilinear Smoluchowski-Poisson system} 

\author{Tomasz Cie\'slak\footnote{Institute of Applied Mathematics, Warsaw University, Banacha 2, 02-097 Warszawa, Poland. \email{T.Cieslak@impan.gov.pl}} \kern8pt \& \kern8pt
Philippe Lauren\c cot\footnote{Institut de 
Math\'ematiques de Toulouse, CNRS UMR~5219, Universit\'e de Toulouse, 118 route de Narbonne, F--31062 Toulouse Cedex 9, 
France. \email{Philippe.Laurencot@math.univ-toulouse.fr}}}
\date{\today}
\maketitle

\begin{abstract}
We prove that, unlike in several space dimensions,  there is no critical (nonlinear) diffusion coefficient for which solutions to the one dimensional quasilinear Smoluchowski-Poisson equation with small mass exist globally while finite time blowup could occur for solutions with large mass.
\end{abstract}

\section{Introduction}\label{i}


In a previous paper \cite{cl3} we investigate the influence of the diffusion coefficient $a$ on the life span of solutions to the one dimensional Smoluchowski-Poisson system
\begin{eqnarray}
\label{he1}
\partial_t{u} &=& \partial_x \left( a(u) \partial_x u -u \partial_x v \right) \;\;\mbox{in}\;\;(0,\infty)\times (0,1),\\
\label{he12}
0&=&\partial_x^2 v + u - M \;\;\mbox{in}\;\;(0,\infty)\times (0,1),\\
\label{he2}
a(u) \partial_x u & = & \partial_x v = 0\;\;\mbox{on} \;\; (0,\infty)\times \{0,1\},\\
\label{he3}
u(0)&=& u_0\ge 0 \;\;\mbox{in} \;\;(0,1),\;\;\int_0^1v(t,x)dx=0\;\;\mbox{for any}\;\; t\in (0,\infty),
\end{eqnarray}
where 
$$
M := \langle u_0 \rangle = \int_0^1 u_0(x) dx
$$ 
denotes the mean value of $u_0$, and uncover a fundamental difference with the quasilinear Smoluchowski-Poisson system in higher space dimensions. More precisely, when the space dimension $n$ is greater or equal to two, there is a critical diffusion $a_*(r):= (1+r)^{(n-2)/2}$ which separates different behaviours for the quasilinear Smoluchowski-Poisson system. Roughly speaking, 
\begin{description}
\item[(a)] if the diffusion coefficient $a$ is stronger than $a_*$ (in the sense that $a(r)\ge C (1+r)^\alpha$ for some $\alpha>(n-2)/n$ and $C>0$), then all solutions exist globally whatever the value of the mass of the initial condition $u_0$ \cite{cw_1},
\item[(b)] if the diffusion coefficient $a$ is weaker than $a_*$ (in the sense that $a(r)\le C (1+r)^\alpha$ for some $\alpha<(n-2)/n$ and $C>0$), then  there exists for all $M>0$ an initial condition $u_0$ with $\langle u_0 \rangle=M$ for which the corresponding solution to the quasilinear Smoluchowski-Poisson system blows up in finite time (in the sense that $\|u(t)\|_\infty\to\infty$ as $t\to T$ for some $T\in (0,\infty)$) \cite{cl1,cw_1,Nagai},
\item[(c)] if the diffusion coefficient $a$ behaves as $a_*$ for large values of $r$, solutions starting from initial data $u_0$ with small mass $\langle u_0 \rangle$ exist globally while there are initial data with large mass for which the corresponding solution to the quasilinear Smoluchowski-Poisson system blows up in finite time \cite{cl1,Nagai}.
\end{description}
Observe that, in space dimension $n=2$, the critical diffusion is constant and a more precise description of the situation {\bf (c)} is actually available. Namely, when $a\equiv 1$, there is a threshold mass $M_*$ such that, if $\langle u_0 \rangle<M_*$, the corresponding solution is global while, for any $M>M_*$, there are initial data with $\langle u_0 \rangle=M$ for which the corresponding solution blows up in finite time \cite{JL92,Nagai,Nagai2}. The threshold mass $M_*$ is known explicitly ($M_*=4\pi$) but it is worth mentioning that for radially symmetric solutions in a ball, the threshold mass is $8\pi$. Similar results are also available for the quasilinear Smoluchowski-Poisson system in $\RR^n$, $n\ge 2$ \cite{B-C-L,ChS04,sugi,sugi2}.

Most surprisingly, the above description fails to be valid in one space dimension and we prove in particular in \cite{cl3} that all solutions are global for the diffusion $a(r)=(1+r)^{-1}$ though it is a natural candidate to be critical. We actually identify two classes of diffusion coefficients $a$ in \cite{cl3}, one for which all solutions exist globally as in {\bf (a)} and the other for which there are solutions blowing-up in finite time starting from initial data with an arbitrary positive mass as in {\bf (b)}, but the situation {\bf (c)} does not seem to occur in one space dimension. The purpose of this note is to show that the dichotomy  {\bf (a)} or  {\bf (b)} can be extended to larger classes of diffusion, thereby extending the analysis performed in \cite{cl3}.

\begin{theo}\label{mainth}
Let the diffusion coefficient $a\in \mathcal{C}^1((0,\infty))$ be a positive function. 
\smallskip

\noindent
(i) Assume first that $a\in L^1(1,\infty)$ and one of the following assumptions is satisfied, either 
\begin{equation}\label{(1)}
\gamma:=\sup_{r\in(0,1)} r\int_r^{\infty}a(s)ds<\infty,
\end{equation}
or there exist $\vartheta>0$ and $\alpha\in (\vartheta/(1+\vartheta),2]$ such that
\begin{equation}\label{decr}
\gamma_\vartheta:=\sup_{r\in (0,1)}{ r^{2+\vartheta} a(r) } < \infty \;\;\mbox{ and }\;\; C_\infty := \sup_{r\ge 1} r^\alpha a(r) < \infty.
\end{equation}
For any $M>0$,  there exists a positive initial condition $u_0\in\mathcal{C}([0,1])$ such that $\langle u_0 \rangle = M$ and the corresponding classical solution to (\ref{he1})-(\ref{he3}) blows up in finite time. 
\smallskip

\noindent
(ii) Assume next that  $a\not\in L^1(1,\infty)$ and consider an initial condition $u_0\in\mathcal{C}([0,1])$ such that $u_0\ge m_0>0$ and $\langle u_0 \rangle = M$ for some $M>0$ and $m_0\in (0,M)$. Then the corresponding classical solution to (\ref{he1})-(\ref{he3}) exists globally. 
\end{theo}

As already mentioned, Theorem~\ref{mainth} extends the results obtained in \cite{cl3}. More precisely, in \cite[Theorem~5]{cl3}, the assertion~(ii) of Theorem~\ref{mainth} is proved under the additional assumption that, for each $\varepsilon\in (0,\infty)$, there is $\kappa_\varepsilon>0$ for which
$$
a(r) \le \varepsilon\ r a(r) + \frac{\kappa_\varepsilon}{r} \;\;\mbox{ for }\;\; r\in (0,1)\,,
$$
which roughly means that $a$ cannot have a singularity stronger than $1/r$ near $r=0$. This assumption turns out to be unnecessary for global existence but nevertheless ensures the global boundedness of the solution in $L^\infty$. Under the sole assumption of Theorem~\ref{mainth}~(ii),  our proof does not exclude that solution to (\ref{he1})-(\ref{he3}) becomes unbounded as $t\to\infty$. Concerning Theorem~\ref{mainth}~(i), it is established in \cite[Theorem~10]{cl3} for $a\in L^1(1,\infty)$ such that there is a concave function $B$ for which
\begin{eqnarray}
0 \le -r A(r) & \le & B(r) \;\;\mbox{ with }\;\; A(r) = - \int_r^\infty a(s)\ ds\,, \qquad r\in (0,\infty)\,, \label{bu1} \\
\lim_{r\to \infty} \frac{B(r)}{r} & = & 0\,. \label{bu2}
\end{eqnarray}
We make this criterion more explicit here by showing that the integrability of $a$ on $(1,\infty)$ and (\ref{(1)}) guarantee the existence of a concave function $B$ satisfying (\ref{bu1}) and (\ref{bu2}), see Lemma~\ref{imp} below. Let us point out here that the assumption (\ref{(1)}) somehow means that $a$ cannot have a singularity stronger that $1/r^2$ near $r=0$. However, the result remains true if $a$ has an algebraic singularity of higher order near $r=0$ which is allowed by (\ref{decr}) provided $a$ decays suitably at infinity. Observe that the second condition in (\ref{decr}) is compatible with the integrability of $a$ at infinity as $\vartheta/(1+\vartheta)<1$. 

Summarizing the outcome of Theorem~\ref{mainth}, we realize that, for a given diffusion coefficient $a$ with a singularity weaker than $1/r^2$ near $r=0$, the integrability or non-integrability of $a$ at infinity completely determines whether we are in the situation {\bf (a)} or {\bf (b)} described above and excludes the situation {\bf (c)}. There is thus no critical diffusion in this class. The same comment applies to the class of diffusion coefficients satisfying (\ref{decr}) with an algebraic singularity stronger than $1/r^2$ near $r=0$. In particular there is no critical nonlinearity in the class of functions ${\cal C}([0,\infty))\cap {\cal C}^1((0,\infty))$.

The paper is organized as follows: in section~\ref{pre} we recall some statements from \cite{cl3}. Section~\ref{ftb} is devoted to proving the finite time blowup of solutions to (\ref{he1})-(\ref{he3}) when $a\in L^1(1,\infty)$. Global existence of solutions for all initial data when $a$ is not integrable at infinity is proved in the last section.

\section{Preliminaries.}\label{pre}

In this section we summarize some results and methods introduced in \cite{cl3}. Let $a\in \mathcal{C}^1((0,\infty))$ be a positive function and consider an initial condition $u_0\in\mathcal{C}([0,1])$ such that $u_0\ge m_0>0$ and $\langle u_0 \rangle = M$ for some $M>0$ and $m_0\in (0,M)$. By \cite[Propositions~2 and~3]{cl3} there is a unique maximal classical solution $(u,v)$ to (\ref{he1})-(\ref{he3}) defined on $[0,T_{max})$ which satisfies
\begin{equation}
\label{zz1}
\min_{x\in [0,1]}{u(t,x)} > 0\,, \quad \langle u(t) \rangle := \int_0^1 u(t,x)\ dx = M\,, \;\;\mbox{ and }\;\; \langle v(t) \rangle := \int_0^1 v(t,x)\ dx = 0
\end{equation}
for $t\in (0,T_{max})$. In addition, $T_{max}=\infty$ or $T_{max}<\infty$ with $\|u(t)\|_\infty\to \infty$ as $t\to T_{max}$.

We next recall the approach introduced in \cite{cl3} which will be used herein as well. Owing to the positivity (\ref{zz1}) and the regularity of $u$, the indefinite integral 
$$
U(t,x) := \int_0^x u(t,z) dz\,, \quad x\in [0,1]\,,
$$
is a smooth increasing function from $[0,1]$ onto $[0,M]$ for each $t\in [0,T_{max})$ and has a smooth inverse $F$ defined by 
\begin{equation}\label{zero}
U(t,F(t,y))=y\,, \qquad (t,y)\in [0,T_{max})\times [0,M]\,.
\end{equation}   
Introducing $f(t,y):=\partial_y F(t,y)$, we have
\begin{equation}\label{odw}
f(t,y)\ u(t,F(t,y))=1\,, \qquad (t,y)\in [0,T_{max})\times [0,M]\,,
\end{equation} 
and it follows from (\ref{he1})-(\ref{he3}) that $f$ solves
\begin{eqnarray}
\partial_t f & = & \partial_y^2 \Psi(f)-1+Mf\,, \qquad (t,y)\in (0,T_{max})\times (0,M)\,, \label{main} \\
\partial_y f(t,0) & = & \partial_y f(t,M)=0\,, \qquad t\in (0,T_{max})\,, \label{boundary} \\
f(0,y) & = & f_0(y) := \frac{1}{u_0(F(0,y))}\,, \qquad y\in (0,M)\,, \label{init}
\end{eqnarray}
where
\begin{equation}\label{Psi}
\Psi'(r):=\frac{1}{r^2}\ a\left( \frac{1}{r} \right) \;\;\mbox{for any}\;\;r>0\,, \qquad \Psi(1):=0\,,
\end{equation} 
Moreover the conservation of mass (\ref{zz1}) yields
\begin{equation}\label{mas}
\int_0^M f(t,y) dy=F(t,M)-F(t,0)=1\,, \qquad t\in [0,T_{max})\,.
\end{equation}
At this point, the crucial observation is that, thanks to (\ref{odw}), finite time blowup of $u$ is equivalent to the vanishing (or touch-down) of $f$. In other words, $u$ exist globally if the minimum of $f(t)$ is positive for each $t>0$. We refer to \cite[Proposition 1]{cl3} for a more detailed description. 

An salient property of  (\ref{he1})-(\ref{he3}) is the existence of a Liapunov function \cite[Lemma~8]{cl3} which we recall now:
\begin{Le}\label{lem}
The function 
$$
L_1(t):= \frac{1}{2} \int_0^M \left| \partial_y\Psi(f(t,y)) \right|^2\ dy+\int_0^M \left( \Psi(f(t,y)) -M\ \Psi_1(f(t,y)) \right)\ dy
$$
is a non-increasing function of time on $[0,T_{max})$, the function $\Psi_1$ being defined by 
\begin{equation}
\label{gex3}
\Psi_1(1):=0 \;\;\mbox{ and }\;\; \Psi_1'(r) := r \Psi'(r) = \frac{1}{r}\ a \left( \frac{1}{r} \right) \,, \qquad r\in (0,\infty)\,.
\end{equation}
\end{Le}

\section{Finite time blowup.}\label{ftb} 

In this section we prove the blowup assertion of Theorem~\ref{mainth}. To this end we first prove that the condition (\ref{(1)}) allows us to construct a concave function $B$ satisfying (\ref{bu1}) and (\ref{bu2}) so that \cite[Theorem~10]{cl3}  can be applied.
\begin{Le}\label{imp}
Let $a\in {\cal C}^1((0,\infty))$ be a positive function such that $a\in L^1(1,\infty)$ and (\ref{(1)}) holds.
Then there exists a concave function $B\in {\cal C}([0,\infty))$ such that for all $r\geq 0$
\begin{equation}\label{(3)}
B(r)\geq r\int_r^{\infty}a(s)ds
\end{equation}
and 
\begin{equation}\label{(4)}
\lim_{r\rightarrow\infty}\frac{B(r)}{r}=0.
\end{equation} 
\end{Le}

\vspace{0.3cm}
\noindent
\textit{Proof of Lemma~\ref{imp}.} We construct $B:[0,\infty)\rightarrow[0,\infty)$ in the following way: we put  
$$
b_i:=\int_{2^i}^\infty a(s)ds\,, \qquad i\ge 0\,,
$$
and notice that $\{b_i\}_{i\ge 0}$ is a decreasing sequence converging to zero as $i\to\infty$. We next define
\begin{equation}
B(r) = \left\{\begin{array}{ccl}
\displaystyle{b_0 r+\gamma}& \mbox{if} &r\in[0,2],\\ 
 & & \\
\displaystyle{b_i r+\sum_{j=0}^{i-1}(b_j-b_{j+1})2^{j+1}+\gamma}& \mbox{if}&r\in(2^i,2^{i+1}] \;\mbox{ and }\; i\geq 1,
\end{array}\right.
\end{equation}
Clearly, $B\in \mathcal{C}([0,\infty))$ and 
\begin{equation}
B'(r) = \left\{\begin{array}{ccl}
b_0 & \mbox{if} &r\in(0,2),\\ 
b_i & \mbox{if} &r\in(2^i,2^{i+1})  \;\mbox{ and }\; i\geq 1.
\end{array}\right.
\end{equation}
Hence $B$ is concave as a consequence of the fact that the sequence $\{b_i\}_{i\ge 0}$ is decreasing. 
Furthermore, for $r\in[0,1]$, we have
\[
B(r)\geq \gamma\geq r\int_r^{\infty}a(s)ds,
\]
and for $r\in[2^i,2^{i+1}]$, $i\geq 0$,
\[
B(r)\geq b_i r=r\int_{2^i}^\infty a(s)ds\geq r\int_r^\infty a(s)ds.
\]
Therefore, $B$ satisfies (\ref{(3)}). 

Finally, let $k\ge 1$. If $i\geq k+1$ and $r\in(2^i,2^{i+1}]$, then
\begin{eqnarray*}
\frac{B(r)}{r} &=& b_i+\frac{\gamma}{r}+\sum_{j=0}^{i-1}(b_j-b_{j+1})\frac{2^{j+1}}{r}\leq b_i+\frac{\gamma}{r}+\sum_{j=k}^{i-1}(b_j-b_{j+1})+\sum_{j=0}^{k-1}(b_j-b_{j+1})\frac{2^{j+1}}{r} \\
&\le & b_i+\frac{1}{r}\left(\gamma+2^k\sum_{j=0}^{k-1}(b_j-b_{j+1})\right)+(b_k-b_i)\leq b_k+\frac{1}{r}\left(\gamma+2^kb_0\right).
\end{eqnarray*}
Consequently,
\[
\limsup_{r\rightarrow\infty}\frac{B(r)}{r} \leq b_k\;\mbox{for all}\;k\geq1\,.
\]
Letting $k\rightarrow\infty$, we obtain (\ref{(4)}) since $b_k\to 0$ as $k\to\infty$ and Lemma~\ref{imp} is proved. \qed

\vspace{0.3cm}
\noindent
\textit{Proof of Theorem~\ref{mainth}~(i), Part~1}. When $a$ belongs to $L^1(1,\infty)$ and satisfies (\ref{(1)}), it follows from Lemma~\ref{imp} that the conditions (\ref{bu1}) and (\ref{bu2}) are satisfied so that Theorem~\ref{mainth}~(i) follows from \cite[Theorem~10]{cl3}.  \qed

\medskip

To handle the other case, we proceed in a different way by showing an upper bound for the function $f$ defined in section~\ref{pre}. We first observe that the function $\Psi$ defined in (\ref{Psi}) satisfies
$$
\Psi(r) = \int_1^\infty \frac{1}{s^2} a\left( \frac{1}{s} \right)\ ds = \int_{1/r}^1 a(s)\ ds\,, \quad r\in (0,\infty)\,,
$$
so that, if $a\in L^1(1,\infty)$, $\Psi(r)$ has a finite limit $\Psi(0):= - \|a\|_{L^1(1,\infty)}$ as $r\to 0$. We then define
\begin{equation}\label{Psit}
\tilde{\Psi}(r) := \Psi(r) - \Psi(0) = \int_0^r \frac{1}{s^2} a\left( \frac{1}{s} \right)\ ds = \int_{1/r}^\infty a(s)\ ds\,, \quad r\in (0,\infty)\,.
\end{equation}

\begin{Le}\label{impo}
Let $a\in {\cal C}^1((0,\infty))$ be a positive function such that $a\in L^1(1,\infty)$. There exists a positive constant $\mu_M>0$ depending only on $M$ and $a$ such that, for any non-negative function $g\in H^1(0,M)$ satisfying $\|g\|_{L^1(0,M)}=1$, we have
\begin{equation}\label{boundb}
\|\tilde{\Psi}(g)\|_{L^\infty(0,M)}^2 \le 32M \mathcal{L}_1(g) + \mu_M,
\end{equation}
with
\begin{equation}\label{pim}
\mathcal{L}_1(g) := \frac{1}{2} \|\partial_y \Psi(g)\|_{L^2(0,M)}^2 + \int_0^M \left( \Psi(g) - M \Psi_1(g) \right)(y)\ dy\,,
\end{equation}
the functions $\Psi$ and $\Psi_1$ being defined in (\ref{Psi}) and (\ref{gex3}), respectively.
\end{Le}

\vspace{0.3cm}
\noindent
\textit{Proof of Lemma~\ref{impo}.} We set $G:= \|g\|_{L^\infty(0,M)}$ which is finite owing to the continuous embedding of $H^1(0,M)$ in $L^\infty(0,M)$. Assume first that $G>1$. Then, for $y\in (0,M)$ and $z\in (0,M)$, we have 
$$
\tilde{\Psi}(g(y)) = \tilde{\Psi}(g(z)) + \int_z^y \partial_{x} \tilde{\Psi}(g(x))\ dx \le \tilde{\Psi}(g(z)) + M^{1/2} \|\partial_y \Psi(g)\|_{L^2(0,M)}.
$$
Integrating the above inequality over $(0,M)$ with respect to $z$ gives
\begin{eqnarray*}
M\tilde{\Psi}(g(y)) & \le & \int_0^M \tilde{\Psi}(g(z))\ dz + M^{3/2} \|\partial_y \Psi(g)\|_{L^2(0,M)} \\
& \le & \int_0^M \mathbf{1}_{[0,2/M]}(g(z)) \tilde{\Psi}(g(z))\ dz + \int_0^M \mathbf{1}_{(2/M,\infty)}(g(z)) \tilde{\Psi}(g(z))\ dz + M^{3/2} \|\partial_y \Psi(g)\|_{L^2(0,M)} \\
& \le & M \tilde{\Psi}\left( \frac{2}{M} \right) + \frac{M \tilde{\Psi}(G)}{2}\ \int_0^M \mathbf{1}_{(2/M,\infty)}(g(z)) g(z)\ dz + M^{3/2} \|\partial_y \Psi(g)\|_{L^2(0,M)} \\
& \le & M \tilde{\Psi}\left( \frac{2}{M} \right) + \frac{M \tilde{\Psi}(G)}{2} + M^{3/2} \|\partial_y \Psi(g)\|_{L^2(0,M)}\,,
\end{eqnarray*}
where we have used the property $\|g\|_{L^1(0,M)}=1$ to obtain the last inequality. Taking the supremum over $y\in (0,M)$ and using the monotonicity and non-negativity of $\tilde{\Psi}$, we deduce that
\begin{equation}\label{y4b}
\tilde{\Psi}(G) \le 2\tilde{\Psi}\left( \frac{2}{M} \right) + 2M^{1/2} \|\partial_y \Psi(g)\|_{L^2(0,M)}.
\end{equation}
We next observe that the integrability of $a$ at infinity also ensures that $\Psi_1(0)>-\infty$, so that $\tilde{\Psi}_1:= \Psi_1-\Psi_1(0)$ is well-defined and satisfies
\begin{equation}\label{bas}
\tilde{\Psi}_1(r) = \int_0^r s \Psi'(s)\ ds \le r \tilde{\Psi}(r)\,, \quad r\in (0,\infty)\,.
\end{equation}
Since $\|g\|_{L^1(0,M)}=1$, it follows from (\ref{y4b}) and (\ref{bas}) that
\begin{equation}\label{y4c}
\int_0^M \tilde{\Psi}_1(g)\ dy \le \int_0^M g \tilde{\Psi}(g)\ dy \le \tilde{\Psi}(G)\ \int_0^M g dy \le 2\tilde{\Psi}\left( \frac{2}{M} \right) + 2M^{1/2} \|\partial_y \Psi(g)\|_{L^2(0,M)}.
\end{equation}
We next infer from (\ref{y4c}) and the non-negativity of $\tilde{\Psi}$ that
\begin{eqnarray*}
\mathcal{L}_1(g) & \ge & \frac{1}{2} \|\partial_y \Psi(g)\|_{L^2(0,M)}^2 + \int_0^M \tilde{\Psi}(g)\ dy + M \Psi(0) - M \int_0^M \tilde{\Psi}_1(g)\ dy \\
& \ge & \frac{1}{2} \|\partial_y \Psi(g)\|_{L^2(0,M)}^2 + M \Psi(0) - 2M \tilde{\Psi}\left( \frac{2}{M} \right) - 2M^{3/2} \|\partial_y \Psi(g)\|_{L^2(0,M)}\\
& \ge & \frac{1}{4} \|\partial_y \Psi(g)\|_{L^2(0,M)}^2 + \left( \frac{1}{2} \|\partial_y \Psi(g)\|_{L^2(0,M)} - 2 M^{3/2} \right)^2 - 4M^3 + M \Psi(0) - 2M \tilde{\Psi}\left( \frac{2}{M} \right) \\
& \ge & \frac{1}{4} \|\partial_y \Psi(g)\|_{L^2(0,M)}^2 - 4M^3 + M \Psi(0) - 2M \tilde{\Psi}\left( \frac{2}{M} \right) ,
\end{eqnarray*}
whence
$$
\|\partial_y \Psi(g)\|_{L^2(0,M)}^2 \le 4 \mathcal{L}_1(g) + 16M^3 - 4M \Psi(0) +8M \tilde{\Psi}\left( \frac{2}{M} \right) .
$$
It then follows from (\ref{y4b}) and the above inequality that 
\begin{eqnarray*}
\tilde{\Psi}(G)^2 & \le & 8 \tilde{\Psi}\left( \frac{2}{M} \right)^2 + 8M \|\partial_y \Psi(g)\|_{L^2(0,M)}^2 \\
& \le & 8 \tilde{\Psi}\left( \frac{2}{M} \right)^2 + 32M \mathcal{L}_1(g) + 128M^4 - 32M^2 \Psi(0) +64M^2 \tilde{\Psi}\left( \frac{2}{M} \right) \\
& \le & 32M \mathcal{L}_1(g) + \mu_M,
\end{eqnarray*}
with 
$$
\mu_M := 1 + 128M^4 - 32M^2 \Psi(0) +64M^2 \tilde{\Psi}\left( \frac{2}{M} \right) + 8 \tilde{\Psi}\left( \frac{2}{M} \right)^2 + \Psi(0)^2 - 32M\Psi(0).
$$
We have thus shown Lemma~\ref{impo} when $G=\|g\|_{L^\infty(0,M)}>1$. To complete the proof, we finally consider the case $G\in [0,1]$ and notice that, in that case, 
$$
0 \le \tilde{\Psi}(G) \le -\Psi(0) \;\;\mbox{ and }\;\; \mathcal{L}_1(g) \ge \int_0^M \tilde{\Psi}(g)\ dy + M \Psi(0) \ge M \Psi(0) ,
$$
since $\Psi_1\le 0$ in $(0,1)$ and $\tilde{\Psi}\ge 0$. Consequently, 
\begin{eqnarray*}
\tilde{\Psi}(G)^2 & \le & \Psi(0)^2 = 32M \Psi(0) + \Psi(0)^2 - 32M\Psi(0) \le 32M \mathcal{L}_1(g) + \mu_M,
\end{eqnarray*} 
and the proof of Lemma~\ref{impo} is complete. \qed

\medskip

As an obvious consequence of Lemmas~\ref{lem} and~\ref{impo} we have the following result:

\begin{Cor}\label{cor:y2}
Let $a\in {\cal C}^1((0,\infty))$ be a positive function such that $a\in L^1(1,\infty)$. For $t\in [0,T_{max})$ and $y\in [0,M]$, we have
$$
0 \le \tilde{\Psi}(f(t,y)) \le \left( 32M \max{\{\mathcal{L}_1(f_0),0\}} + \mu_M \right)^{1/2}.
$$
\end{Cor}

\vspace{0.3cm}
\noindent
\textit{Proof of Corollary~\ref{cor:y2}}. Clearly $\mathcal{L}_1(f(t))=L_1(t)\le L_1(0)=\mathcal{L}_1(f_0)\le \max{\{\mathcal{L}_1(f_0),0\}}$ for $t\in [0,T_{max})$ by Lemma~\ref{lem} and Corollary~\ref{cor:y2} readily follows from Lemma~\ref{impo}. \qed

\medskip

\begin{Rem}\label{rem:positivity}
Corollary~\ref{cor:y2} provides an $L^\infty$-bound on $f$ only if $\Psi(r)\to\infty$ as $r\to\infty$, that is, if $a\not\in L^1(0,1)$. In that case, it gives a positive lower bound for $u$ by (\ref{odw}).
\end{Rem}

We next turn to the proof of the second part of Theorem~\ref{mainth} for which we develop further the arguments from \cite[Theorem~10]{cl3}. 

\vspace{0.3cm}
\noindent
\textit{Proof of Theorem~\ref{mainth}~(i), Part~2}. Assume now that $a\in L^1(1,\infty)$ and satisfies (\ref{decr}). We fix $M>0$, $q>2$, and $\varepsilon_M\in (0,1)$ such that 
\begin{equation}\label{prunelle}
q> \max{\left\{ 3+\vartheta, \frac{5+3\vartheta}{\alpha(\vartheta+1)-\vartheta} \right\}} \;\;\;\mbox{ and }\;\;\;  \frac{q(q+1)}{M^2} \int_{1/\varepsilon_M}^\infty a(s)\ ds \le \frac{1}{2}\,,
\end{equation} 
the existence of $\varepsilon_M$ being guaranteed by the integrability of $a$ at infinity. 

For 
\begin{equation}\label{gaston}
\delta\in \left( 0,\min{\left\{ 1 , 2M ,(2M)^{-1/q} \right\}} \right)\,,
\end{equation}
we put
\begin{equation}\label{pam}
f_0(y) := \frac{2(1-M\delta^q)}{\delta^2}\ (\delta-y)_+ + \delta^q\ge \delta^q>0\,, \quad y\in [0,M]\,.
\end{equation}
Then
\begin{equation}\label{poum}
\int_0^M f_0(y)\ dy = 1, \quad \|f_0\|_{L^\infty(0,M)} = \frac{2(1-M\delta^q)}{\delta} + \delta^q\le \frac{2}{\delta}\,.
\end{equation}
Introducing next 
$$
m_q(t) := \int_0^M y^q f(t,y)\ dy\,, \quad t\in [0,T_{max})\,,
$$
we have
\begin{equation}\label{boule}
m_q(0) = \left( \frac{2(1-M\delta^q)}{(q+1)(q+2)} + \frac{M^{q+1}}{q+1}\right)\ \delta^q \le C_1 \delta^q \;\;\mbox{ with }\;\; C_1 := \left( \frac{2 +(q+2) M^{q+1}}{(q+1)(q+2)} \right)\,.
\end{equation}

It follows from (\ref{main}), (\ref{boundary}), and the non-negativity of $\tilde{\Psi}$ that 
\begin{eqnarray}
\frac{dm_q}{dt} & = & - q\ \int_0^M y^{q-1} \partial_y \tilde{\Psi}(f)\ dy + M m_q - \frac{M^{q+1}}{q+1}, \nonumber \\
\frac{dm_q}{dt} & \le & q(q-1)\ \int_0^M y^{q-2} \tilde{\Psi}(f)\ dy + M m_q - \frac{M^{q+1}}{q+1} \label{y12}\,.
\end{eqnarray}
We shall now estimate the integral on the right-hand side of (\ref{y12}): to this end, we split the domain of integration into three parts which we handle differently. As a preliminary step, we notice that, by (\ref{decr}), 
\begin{equation}\label{y12b}
\Psi'(r) \le \gamma_\vartheta r^\vartheta \;\;\;\mbox{ and }\;\;\; \Psi(r) \le \frac{\gamma_\vartheta}{\vartheta+1}\ r^{\vartheta+1} \le \gamma_\vartheta\ r^{\vartheta+1}\,, \quad r\ge 1\,.
\end{equation}
We next define 
$$
K_0:= \left( 32M \max{\{\mathcal{L}_1(f_0),0\}} + \mu_M \right)^{1/2(2+\vartheta)} > 1\,,
$$
and consider $(t,y)\in [0,T_{max})\times [0,M]$. 
\begin{itemize}
\item If $f(t,y)\in (0,\varepsilon_M]$, it follows from (\ref{prunelle}) and the monotonicity of $\tilde{\Psi}$ that  
\begin{equation}\label{jeden}
\tilde{\Psi}(f(t,y)) \le \tilde{\Psi}(\varepsilon_M) = \int_{1/\varepsilon_M}^\infty a(s)\ ds \le  \frac{M^2}{2q(q+1)}\,.
\end{equation}
\item If $f(t,y)\in (\varepsilon_M,K_0)$, then (\ref{y12b}) and the monotonicity of $\tilde{\Psi}$ yield
\begin{equation}\label{dwa}
\tilde{\Psi}(f(t,y)) = \frac{\tilde{\Psi}(f(t,y))}{f(t,y)}\ f(t,y) \le  \frac{\Psi(K_0)-\Psi(0)}{\varepsilon_M}\ f(t,y) \le \frac{\gamma_\vartheta K_0^{\vartheta+1}  - \Psi(0)}{\varepsilon_M}\ f(t,y)\,.
\end{equation}
\item If $f(t,y)\ge K_0$, Corollary~\ref{cor:y2} ensures that 
\begin{equation}\label{trzy}
\tilde{\Psi}(f(t,y)) = \frac{\tilde{\Psi}(f(t,y))}{f(t,y)}\ f(t,y) \le  \frac{K_0^{\vartheta+2}}{K_0}\ f(t,y) \le K_0^{\vartheta+1}\ f(t,y)\,.
\end{equation}
\end{itemize}
Consequently, recalling that $K_0>1$ and $\Psi(0)<0$, we deduce from (\ref{y12}) and (\ref{jeden})-(\ref{trzy}) that 
\begin{eqnarray*}
\frac{dm_q}{dt} & \le &  q(q-1)\ \int_0^M y^{q-2} \tilde{\Psi}(f)\ \mathbf{1}_{(0,\varepsilon_M]}(f)\ dy + q(q-1)\ \int_0^M y^{q-2} \tilde{\Psi}(f)\ \mathbf{1}_{(\varepsilon_M,K_0)}(f)\ dy \\
& + & q(q-1)\ \int_0^M y^{q-2} \tilde{\Psi}(f)\ \mathbf{1}_{[K_0,\infty)}(f)\ dy + M m_q - \frac{M^{q+1}}{q+1} \\
& \le & \frac{(q-1)M^2}{2(q+1)}\ \int_0^M y^{q-2}\ dy + q(q-1)\ \frac{\gamma_\vartheta K_0^{\vartheta+1}  - \Psi(0)}{\varepsilon_M}\ \int_0^M y^{q-2} f\ dy \\
& + & q(q-1) K_0^{\vartheta+1}\ \int_0^M y^{q-2} f\ dy + M m_q - \frac{M^{q+1}}{q+1} \\
& \le & C_2\  K_0^{\vartheta+1}\ \int_0^M y^{q-2} f\ dy + M m_q - \frac{M^{q+1}}{2(q+1)}\,,
\end{eqnarray*}
with $C_2:= q(q-1)(\gamma_\vartheta-\Psi(0)+\varepsilon_M)/\varepsilon_M$. We next use H\"older's inequality and (\ref{mas}) to conclude that
\begin{equation}\label{jules}
\frac{dm_q}{dt} \le C_2\ K_0^{\vartheta+1}\ m_q^{(q-2)/q} + M m_q - \frac{M^{q+1}}{2(q+1)}\,.
\end{equation}

It remains to estimate $K_0$ and in fact $\mathcal{L}_1(f_0)$. Since $\Psi$ is negative on $(0,1)$ and $\Psi_1$ is bounded from below by $\Psi_1(0)$, it follows from (\ref{gaston}) and (\ref{pam}) that 
\begin{eqnarray*}
\mathcal{L}_1(f_0) & \le & \frac{2}{\delta^4}\ (1-M\delta^q)^2\ \int_0^\delta |\Psi'(f_0)|^2\ dy + \int_0^\delta \Psi(f_0)\ dy - M^2 \Psi_1(0) \\
& \le & \frac{2}{\delta^4}\ \int_0^\delta |\Psi'(f_0)|^2\ dy + \int_0^\delta \Psi(f_0)\ dy - M^2 \Psi_1(0) .
\end{eqnarray*}
On the one hand, we infer from (\ref{poum}), (\ref{y12b}), and the monotonicity of $\Psi$ that
$$
\int_0^\delta \Psi(f_0)\ dy \le \delta\ \Psi\left( \frac{2}{\delta} \right) \le \gamma_\vartheta 2^{\vartheta+1}\ \delta^{-\vartheta}\,.
$$
On the other hand, we have
\begin{eqnarray*}
f_0(y) & \ge & 1 \;\;\mbox{ for }\;\; y\in [0,y_\delta] \;\;\mbox{ with }\;\; y_\delta := \delta - \frac{1-\delta^q}{2(1-M\delta^q)} \delta^2>0\,, \\
f_0(y) & \in & [\delta^q,1] \;\;\mbox{ for }\;\; y\in [y_\delta,\delta]\,,
\end{eqnarray*}
so that, if $y\in [0,y_\delta]$,
$$
\Psi'(f_0(y))^2 \le \gamma_\vartheta^2 f_0(y)^{2\vartheta} \le \gamma_\vartheta 4^\vartheta\ \delta^{-2\vartheta}
$$
by (\ref{poum}) and (\ref{y12b}), while, if $y\in (y_\delta,\delta]$,
$$
\Psi'(f_0(y))^2 \le \frac{1}{f_0(y)^4} a\left( \frac{1}{f_0(y)} \right)^2 \le C_\infty^2\ f_0(y)^{2(\alpha-2)} \le C_\infty^2\ \delta^{-2q(2-\alpha)}
$$
by (\ref{decr}) since $\alpha\le 2$. Therefore,
\begin{eqnarray*}
\mathcal{L}_1(f_0) & \le & \frac{2}{\delta^4}\ \left[ \int_0^{y_\delta} \gamma_\vartheta 4^\vartheta\ \delta^{-2\vartheta}\ dy + \int_{y_\delta}^\delta C_\infty^2\ \delta^{-2q(2-\alpha)}\ dy \right] + \gamma_\vartheta 2^{\vartheta+1}\ \delta^{-\vartheta} - M^2 \Psi_1(0)\\
& \le & \gamma_\vartheta 4^{\vartheta+1}\ \delta^{-3-2\vartheta} + C_\infty^2\ \frac{1-\delta^q}{2(1-M\delta^q)}\ \delta^{-2-2q(2-\alpha)} + \gamma_\vartheta 2^{\vartheta+1}\ \delta^{-\vartheta} - M^2 \Psi_1(0)\\
& \le & \gamma_\vartheta 4^{\vartheta+1}\ \delta^{-2(2+\vartheta)} + C_\infty^2\ \delta^{-2-2q(2-\alpha)} + \gamma_\vartheta 2^{\vartheta+1}\ \delta^{-\vartheta} - M^2 \Psi_1(0) \\
& \le & C_3\ \left( \delta^{-2(2+\vartheta)} + \delta^{-2-2q(2-\alpha)}  \right)
\end{eqnarray*}
with $C_3:= \gamma_\vartheta 4^{\vartheta+2} + C_\infty^2 - M^2 \Psi_1(0)$. Therefore,
\begin{equation}\label{cubitus}
K_0^{\vartheta+1} \le C_4\ \left( \delta^{-(\vartheta+1)} + \delta^{-(\vartheta+1)(1+q(2-\alpha))/(\vartheta+2)} \right)
\end{equation}
for some constant $C_4>0$ depending only on $M$ and $a$. 

Combining (\ref{jules}) and (\ref{cubitus}) yields
\begin{equation}\label{senechal}
\frac{dm_q}{dt} \le \Lambda_\delta(m_q) := C_5\ \left( \delta^{-(\vartheta+1)} + \delta^{-(\vartheta+1)(1+q(2-\alpha))/(\vartheta+2)} \right)\ m_q^{(q-2)/q}\ + M m_q - \frac{M^{q+1}}{2(q+1)}
\end{equation}
for $t\in [0,T_{max})$ and some constant $C_5>0$ depending only on $M$ and $a$. At this point, we note that the monotonicity of $\Lambda_\delta$ and (\ref{senechal}) imply that $\Lambda_\delta(m_q(t))\le \Lambda_\delta(m_q(0))$ for $t\in [0,T_{max})$ if $\Lambda_\delta(m_q(0))<0$, the latter condition being satisfied for $\delta$ small enough as 
$$
\Lambda_\delta(m_q(0)) \le C_1^{(q-2)/q} C_5\ \left( \delta^{q-3-\vartheta} + \delta^{(q(\alpha(\vartheta+1)-\vartheta)-3\vartheta-5)/(\vartheta+2)} \right) + M C_1\ \delta^q - \frac{M^{q+1}}{2(q+1)}
$$
by (\ref{prunelle}) and (\ref{boule}).

Summarizing, we have shown that, if $\delta$ satisfies (\ref{gaston}) and
\begin{equation}\label{semaphore}
C_1^{(q-2)/q} C_5\ \left( \delta^{q-3-\vartheta} + \delta^{(q(\alpha(\vartheta+1)-\vartheta)-3\vartheta-5)/(\vartheta+2)} \right) + M C_1\ \delta^q < \frac{M^{q+1}}{2(q+1)}\,,
\end{equation}
we have 
$$
\frac{dm_q}{dt}(t) \le \Lambda_\delta(m_q(t)) \le \Lambda_\delta(m_q(0))<0\,, \quad t\in [0,T_{max})\,,
$$
an inequality which can only be true on a finite time interval owing to the non-negativity of $m_q$. Therefore, $T_{max}<\infty$ in that case and, for any $M>0$, we have found an initial condition $u_0$ given by (\ref{zero}), (\ref{odw}), and (\ref{pam}) (for $\delta$ small enough according to the above analysis) such that $\langle u_0\rangle=M$ and the first component $u$ of the corresponding solution to (\ref{he1})-(\ref{he3}) blows up in finite time. \qed

\section{Global existence.}\label{ge}

The proof of Theorem~\ref{mainth}~(ii) also relies on the study of the function $L_1$ defined in Lemma~\ref{lem}. For that purpose, we first recall another property  from \cite{cl3}. We define the function $E_1$ by 
\begin{equation}\label{landau}
E_1(h) := \frac{1}{2} \|\partial_y h\|_2^2 + \int_0^M \mathbf{1}_{(-\infty,0)}(h(y))\ h(y)\ dy\,, \quad h\in H^1(0,M)\,,
\end{equation}
for which we have the following lower bound.
\begin{Le}\cite[Lemma 9]{cl3}\label{lab1}
For $M>0$, we have
\begin{equation}\label{gex5}
E_1(h) \ge \frac{1}{4}\ \|\partial_y h\|_2^2 - M^3 - M \left| \Psi\left( \frac{1}{M} \right) \right|\,,
\end{equation}
and
\begin{equation}\label{gex6}
\|h\|_1\leq M^{3/2} \|\partial_y h\|_2 + M \left| \Psi\left( \frac{1}{M} \right) \right|
\end{equation} 
for every $h\in H^1(0,M)$ satisfying
\begin{equation}
\label{gex7}
\int_0^M \Psi^{-1}(h)(y)\ dy = 1\,.
\end{equation}
\end{Le}
  
We now show that the non-integrability of $a$ at infinity allows us to show that $T_{max}=\infty$. To this end, we use the alternative formulation (\ref{main})-(\ref{init}) as in \cite{cl3} and prove that $f$ cannot vanish in finite time.

\vspace{0.3cm}
\noindent
\textit{Proof of Theorem~\ref{mainth}~(ii).}  Owing to (\ref{init}) and the assumptions made on $u_0$, we have
\[
0<f_0(y)\leq \frac{1}{m_0}\,, \quad y\in [0,M].
\]
Introducing $\Sigma(t):=M^{-1}+e^{Mt}\left(m_0^{-1}-M^{-1}\right)$ for $t\geq 0$ we have
\[
\partial_t\Sigma-\partial_y^2\Psi(\Sigma)-M\Sigma+1=M\left(\Sigma- \frac{1}{M}\right)-M\Sigma+1=0,
\]
\[
\Sigma(0)=\frac{1}{m_0}\geq f_0(y), \quad y\in(0,M),
\]
and the comparison principle warrants that
\begin{equation}\label{por}
f(t,y)\leq\Sigma(t),\quad (t,y)\in [0,T_{max})\times[0,M].
\end{equation}

We now follow the strategy of the proof of \cite[Theorem 5]{cl3} and first use the properties of $\Psi$, $\Psi_1$, and (\ref{por}) to estimate the function $L_1$ defined in Lemma~\ref{lem} from below. Indeed, since $\Psi\geq 0$ on $(1,\infty)$ and $\Psi_1\leq 0$ on $(0,1)$ we arrive at
\begin{eqnarray*}
L_1(0)\ge L_1(t) &=& \frac{1}{2}\left\|\partial_y\Psi(f(t))\right\|_{L^2(0,M)}^2+\int_0^M \mathbf{1}_{(0,1)}(f(t,y))(\Psi-M\Psi_1)(f(t,y))\ dy \\
& + & \int_0^M \mathbf{1}_{(1,\infty)}(f(t,y))(\Psi-M\Psi_1)(f(t,y))\ dy \\
& \geq & \frac{1}{2}\left\|\partial_y\Psi(f(t))\right\|_{L^2(0,M)}^2+\int_0^M{\bf 1}_{(-\infty,0)}(\Psi(f(t,y)))\Psi(f(t,y))\ dy \\
& - & M\int_0^M{\bf 1}_{(1,\infty)}(f(t,y))\Psi_1(f(t,y))\ dy \\
& \geq & E_1(\Psi(t))-M^2\Psi_1(\Sigma(t)),
\end{eqnarray*}
where $E_1$ is defined in (\ref{landau}) and we have used (\ref{por}) to obtain the last inequality. Next, by Lemma~\ref{lab1} and (\ref{mas}), we have 
$$
L_1(0)\geq \frac{1}{4}\left\|\partial_y\Psi(f(t))\right\|_{L^2(0,M)}^2-M^3-M\left|\Psi\left(\frac{1}{M}\right)\right|-M^2\Psi_1(\Sigma(t)),
$$
whence
\begin{equation}\label{prandtl}
\frac{1}{4}\left\|\partial_y\Psi(f(t))\right\|_{L^2(0,M)}^2\leq L_1(0)+M^3+M\left|\Psi\left(\frac{1}{M}\right)\right|+M^2\Psi_1(\Sigma(t)).
\end{equation}
Using again Lemma~\ref{lab1}, we have
\begin{eqnarray*}
\left\|\Psi(f(t))\right\|_{L^1(0,M)} & \leq & M^{3/2}\left\|\partial_y\Psi(f(t))\right\|_{L^2(0,M)}+M\left|\Psi\left(\frac{1}{M}\right)\right| \\
& \leq & 2M^{3/2}\left(L_1(0)+M^3+M\left|\Psi\left(\frac{1}{M}\right)\right|+M^2\Psi_1(\Sigma(t))\right)^{1/2}+M\left|\Psi\left(\frac{1}{M}\right)\right|.
\end{eqnarray*}
Combining the previous inequality with (\ref{prandtl}) and the Poincar\'e inequality leads us to the bound
\begin{equation}\label{last}
\left\|\Psi(f(t))\right\|_{H^1(0,M)}\leq C_6(T)\,, \quad t\in [0,T]\cap [0,T_{max})\,,
\end{equation}
for all $T>0$. Together with the continuous embedding of $H^1(0,M)$ in $L^\infty(0,M)$, (\ref{last}) gives
\[
-C_7(T)\leq\Psi(f(t,y))\leq C_7(T)\,, \quad (t,y)\in ([0,T]\cap [0,T_{max})) \times [0,M]\,.
\]
Since
$$
\lim_{r\rightarrow 0}\Psi(r)=-\infty
$$ 
due to $a\not\in L^1(1,\infty)$, the above lower bound on $\Psi(f)$ ensures that $f(t)$ cannot vanish in finite time, from which Theorem~\ref{mainth}~(ii) follows as already discussed in section~\ref{pre}. \qed

\noindent\textbf{Acknowledgements.} This work was done while the first author held a post-doctoral position at the University of Z\"{u}rich and is also partially supported by the Polish Ministry of Science and Higher Education under grant number NN201 396937 (2009-2012).


\end{document}